# ETHOS.FINE: A Framework for Integrated Energy System Assessment

**Theresa Groß[1], Kevin Knosala[1,2], Maximilian Hoffmann[1], Noah Pflugradt[1], and Detlef Stolten[1,2]**

**1** Forschungszentrum Jülich GmbH, Institute of Energy and Climate Research – Techno-economic Systems Analysis (IEK-3), 52425 Jülich, Germany
**2** RWTH Aachen University, Chair for Fuel Cells, Faculty of Mechanical Engineering, 52062 Aachen, Germany

## Summary

The decarbonization of energy systems worldwide requires a transformation of their design and operation across all sectors, that is, the residential and commercial, industrial, and transportation sectors. Energy system models are frequently employed for assessing these changes, providing scenarios on potential future system design and on how new technologies and a modified infrastructure will meet future energy demands. Thus, they support investment decisions and policy-making. The Python-based Framework for Integrated Energy System Assessment (ETHOS.FINE) is a software package that provides a toolbox for modeling, analyting and evaluating such energy systems using mathematical optimization. ETHOS.FINE is part of the Energy Transformation paTHway Optimization Suite (ETHOS)[1], a collection of modeling tools developed by the Institute of Energy and Climate Research - Techno-Economic System Analysis (IEK-3) at Forschungszentrum Jülich. ETHOS offers a holistic view on energy systems at arbitrary scales providing tools for geospatial renewable potential analyses, time series simulation tools for residential and industrial sector, discrete choice models for the transportation sector, modeling of global energy supply routes, and local infrastructure assessments, among others. The ETHOS model suite is, e.g., used for analyzing the energy transition of Germany (Stolten et al., 2022).

ETHOS.FINE is not limited to a single instance of energy system. Instead, it can be freely adapted to consider multiple commodities, regions, time steps and investment periods. The optimization objective is to minimize the total annual cost of the system and is subject to technical and environmental constraints. The generic object-oriented implementation allows for arbitrary spatial scales and numbers of regions – from the local level, e.g., individual buildings, to the regional one, e.g., districts or industrial sites, to the national and international levels. Furthermore, the spatial technological resolution can be aggregated by built-in aggregation methods that are described in Patil et al. (2022) for reducing model size and complexity.

This also applies to the model's temporal resolution. Apart from using the full temporal resolution defined by the input data, integrated time series aggregation methods utilizing the built-in Python package tsam[2] allow for reducing the model's complexity and its computation time (M. Hoffmann et al., 2022) while still enabling the flexibility of seasonal storage technologies, despite the reduced model complexity (Kotzur et al., 2018). In addition, ETHOS.FINE enables the investigation of transformation pathways by considering multiple investment periods in a perfect foresight approach, as well as the stochastic optimization for a single year optimization with multiple sets of input parameters (e.g., changing energy demand forecasts or weather conditions) to find more robust energy system designs.

---

[1] ETHOS – Energy Transformation paTHway Optimization Suite, https://www.fz-juelich.de/en/iek/iek-3/expertise/model-services
[2] tsam – Time Series Aggregation Module, https://github.com/FZJ-IEK3-VSA/tsam

## Methodology

ETHOS.FINE comprises seven main classes: the EnergySystemModel class can be seen as the model's container, collecting all relevant input data for its setup. All technologies to be considered are added to this container. The Component class contains the variables and constraints common to all system components, e.g., capacity limits and limitations to the operation of technologies. The five classes of Source, Sink, Conversion, Transmission, and Storage provide the functionalities to model energy generation and consumption, conversion processes, energy storage for later usage, and energy transfer between regions. Each class introduces a specific set of constraints that is added to the optimization program. Supplementary subclasses offer additional component features, e.g., the option to model partial load behavior and ramping constraints for power plants. The described structure is depicted in *Figure 1*.

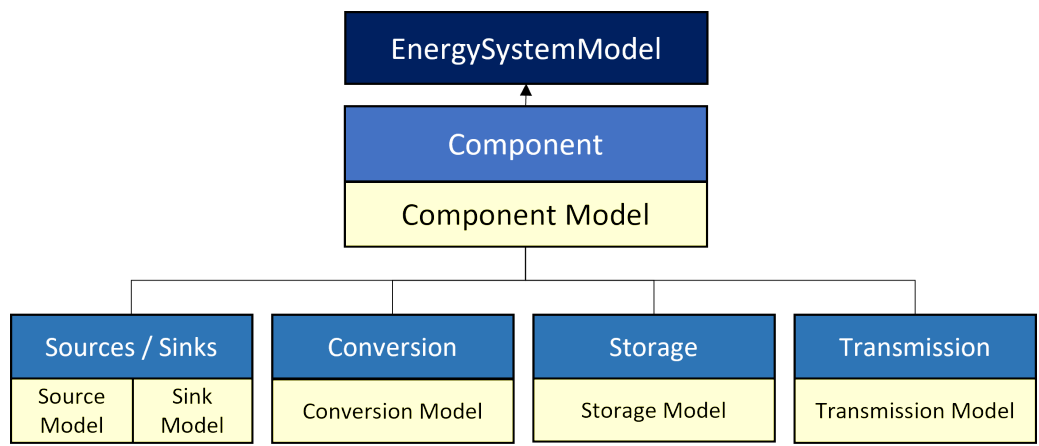

Figure 1 Structure of the main classes in ETHOS.FINE. Additional model classes contain the definition of the specific variables, sets and constraints for each class to set up the optimization model.

The energy system model can be set up as a linear program (LP), or as a mixed integer linear program (MILP), depending on the chosen representation of the added components. The optimization program is a Pyomo[3] instance to enable a flexible solver choice, i.e., ETHOS.FINE optimize energy systems using both, commercial and open source solvers. A description of the variables, constraints and the objective function is given in Welder (2022). Depending on the spatial and temporal resolution of the modeled system, the input parameters are primarily given as Pandas.DataFrames[4] with regions and time steps serving as indices and columns. The model output yields detailed information on the required investments in each region for the installation and operation of the chosen components, as well as the temporally- resolved operation of every component. This also includes the charging and discharging of storage components and commodity flows between regions via transmission components. In addition, the framework also provides plotting options for spatially- and temporally- resolved results. Model input and output can be saved to netCDF files to support reproducibility.

---

[3] Pyomo, Pyoton Optimization Modeling Language, https://pyomo.org/ (Bynum et al., 2021; Hart et al., 2011)

[4] Pandas, Python Data Analysis Library, https://pandas.pydata.org/,

# Statement of need

Groissböck (2019) presents an overview of open source energy system optimization tools and lists their respective integrated functionalities. The author demonstrates their competitiveness against commercial software and outlines hitherto unimplemented functionalities of the described tools. ETHOS.FINE provides a unique generic model setup with a high level of freedom for model developers. Beyond energy system models, its generic implementation enables the modeling of all kinds of optimization issues, such as material flows and resource consumption or conversion as part of life-cycle analysis. The software exhibits many of the features described by Groissböck (2019) and is under constant development with five major releases in the last five years. Its code is openly accessible on Github which allows for contributions and feedback from a wider modeling community. The use cases described in the next section demonstrate the broad range of analyses that can be conducted with the tool. ETHOS.FINE is designed to be used by researchers, students, and for teaching purposes in the field of energy systems modeling. In particular, its exceptional capabilities with respect to complexity reduction (Kotzur et al., 2021) using spatial (Patil et al., 2022) and temporal aggregation (M. Hoffmann et al., 2020, 2021, 2022; M. A. C. Hoffmann, 2023), as well as heuristics to deal with MILPs (Kannengießer et al., 2019; Singh et al., 2022) open up a wide field of applications from small- to global-scale energy system models. For newcomers who are not familiar with programming, it also has the flexibility to set up models by means of Excel files.

# Examples for previous usage

ETHOS.FINE has been used in various studies for energy systems analyses on different scales, leveraging its capability to dynamically adapt to computational complexity. Initial applications can be found in Welder et al. (2018) and Welder et al. (2019): The authors analyzed hydrogen electricity reconversion pathways in a multi-regional energy system model implemented in ETHOS.FINE for the northern part of Germany. Groß (2023) later used the framework to model the future energy system of Germany with a high spatial resolution and thereby to investigate the need for new infrastructure. D. G. Caglayan et al. (2019) set up an ETHOS.FINE model of the European energy system, and analyzed the influence of varying weather years on the cost-optimal system design based on 100% renewable energy source usage. Their findings are also used to determine a robust variable renewable energy sources-based system design, ensuring supply security for a wide range of weather years (D. Caglayan et al., 2021). Knosala et al. (2021) assessed hydrogen technologies in residential buildings in a multi-commodity single-building model. The building model from this work was also used for a sensitivity analysis of energy carrier costs for the application of hydrogen in residential buildings (Knosala et al., 2022). Spiller et al. (2022) analyzed carbon emission reduction potentials for hotels on energetically self-sufficient islands. More recently, Weinand et al. (2023) used the framework to assess the Rhine Rift Valley with respect to its potential for lithium extraction from deep geothermal plants. Meanwhile, Jacob et al. (2023) investigated the potential of Carnot batteries in the German electricity system. Busch et al. (2023) analyzed the role of liquid hydrogen, likewise on a national scale, whereas Franzmann et al. (2023) looked at the green hydrogen cost potentials for global trade. These examples illustrate the manifold application cases that can be tackled by ETHOS.FINE.

## Acknowledgements

We acknowledge contributions from Lara Welder, Robin Beer, Johannes Behrens, Julian Belina, Toni Busch, Arne Burdack, Henrik Büsing, Dilara Caglayan, Philipp Dunkel, David Franzmann, Patrick Freitag, Heidi Heinrichs, Shitab Ishmam, Timo Kannengießer, Leander Kotzur, Stefan Kraus, Felix Kullmann, Dane Lacey, Jochen Linssen, Lars Nolting, Rachel Maier, Shruthi Patil, Jan Priesmann, Oliver Rehberg, Stanley Risch, Martin Robinius, Thomas Schöb, Bismark Singh, Andreas Smolenko, Maximilian Stargardt, Peter Stenzel, Chloi Syranidou, Johannes Thürauf, Christoph Winkler and Michael Zier during the development of this software package.

This work was initially supported by the Helmholtz Association under the Joint Initiative "Energy System 2050 - A Contribution of the Research Field Energy." The authors also gratefully acknowledge financial support by the Federal Ministry for Economic Affairs and Energy of Germany as part of the project METIS (project number 03ET4064, 2018-2022).

This work was supported by the Helmholtz Association under the program Energy System Design.